\documentclass[11pt]{article}
\usepackage{}
\usepackage{amsfonts}
\usepackage{amssymb}
\usepackage{amsmath,amssymb,graphicx}
\usepackage{epsfig}
\newtheorem{theorem}{Theorem}[section]
\newtheorem{theorema}{Theorem}

\newtheorem{lemma}[theorem]{Lemma}
\newtheorem{lemmaa}[theorema]{Lemma}

\newtheorem{corollary}[theorem]{Corollary}

\newcommand{\ma}{\boldsymbol{a}}
\newcommand{\mj}{\boldsymbol{j}}
\newcommand{\me}{\boldsymbol{e}}
\newcommand{\mtheta}{\boldsymbol{\theta}}
\newcommand{\mphi}{\boldsymbol{\phi}}
\newcommand{\kz}{\kappa(z)}
\newcommand{\nn}{\nonumber}

\newcommand{\be}{\begin{eqnarray}}
\newcommand{\ee}{\end{eqnarray}}
\newcommand{\ba}{\begin{array}}
\newcommand{\ea}{\end{array}}
\newcommand{\ben}{\begin{eqnarray*}}
\newcommand{\een}{\end{eqnarray*}}

\newcommand{\pa}{\partial}

\begin{document}

\hoffset -34pt

\title{{\large \bf A note on the asymptotic behavior of conformal metrics with negative
curvatures near isolated singularities}
% \thanks{  }
\thanks{\mbox{Keywords. Singularities, curvatures, metrics.}}}
\author{\normalsize  Tanran Zhang }
\date{}
\maketitle \baselineskip 21pt
\noindent

\begin{minipage}{138mm}
\renewcommand{\baselinestretch}{1} \normalsize
\begin{abstract}
{The asymptotic behavior of conformal metrics with negative curvatures near an isolated singularity for at most second order derivatives was described by Kraus and Roth in one of their papers in 2008. Our work improves one estimate of theirs and shows the estimate for higher order derivatives near an isolated singularity by means of potential theory. We also give some limits of Minda-type for SK-metrics near the origin. Combining these limits with the Ahlfors' lemma, we provide some observations SK-metrics.}
\end{abstract}
\end{minipage}\\
\\
\\
\renewcommand{\baselinestretch}{1} \normalsize
\section{Introduction}

\vspace*{4mm}
The research of conformal metrics has a long history, since the time of Liouville and Picard, see \cite{Liouville1, Picard,Picard2}. For a conformal metric $\lambda(z)|dz|$ on a subdomain $G$ of the complex plane $\mathbb{C}$, we can define its (generalized) Gaussian curvature $\kappa_{\lambda}(z)$. Let $u(z)= \log \lambda(z)$. If $\kappa_{\lambda}(z)=0$, then $u(z)$ satisfies the Laplace equation $\Delta u=0$, which means $u(z)$ is harmonic on $G$. So that the property of $u(z)$ can be studied by means of potential theory, see, e.g. \cite{PDE}. If  $\kappa_{\lambda}(z)=-4$, then $\Delta \log \lambda= 4\lambda^2$ and $u(z)$ is the solution to the Liouville equation
\be\label{liouville equation}
\Delta u=4e^{2u}.
\ee
Each solution to equation (\ref{liouville equation}) belongs to a class of subharmonic functions and it is corresponding to a kind of special metric, called the SK-metric, according to Heins, see \cite{Heins}. The existence and the uniqueness of the solutions to equation (\ref{liouville equation}) are subject to some suitable boundary conditions. Through out our study, we are concerned only with the asymptotic behavior near an isolated singularity of the solution to equation (\ref{liouville equation}), so it is sufficient to consider the behavior in the punctured unit disk $\mathbb{D} \backslash\{0\}$, where the origin is an isolated singularity of some order $\alpha\leq 1$. Near the singularity, we need some more refined invariants to estimate the asymptotic behavior, like the growth of the density. The assignment of the order of the singularity is such an invariant.
\par As for equation (\ref{liouville equation}), Liouville proved that, in any disk $D$ contained in the punctured unit disk $\mathbb{D}\backslash\{0\}$ every solution $u$ to \eqref{liouville equation} can be written as $$u(z)=\log\frac{|f'(z)|}{1-|f(z)|^2},$$
where $f$ is a holomorphic function in $D$, see \cite{Roth1}. Based on  Liouville's results, Nitsche described the behavior of $u(z)$ with constant curvature $\kappa(z)\equiv -4$ near the isolated singularities on plane domains in \cite{Nitsche}. Subsequently, Kraus and Roth extended Nitsche's results to the solutions of the more general equation
\be\label{general equation}
\Delta u=-\kz e^{2u}
\ee
with strictly negative, H\"{o}lder continuous curvature functions $\kz$ in \cite{Rothbehaviour}. In fact, equation (1.2) has an exquisite geometric interpretation: Every solution $u$ to \eqref{general equation} induces a conformal metric $e^{u(z)}|dz|$ with Gaussian curvature function $\kz$ and vice versa (for more details, see \cite{Rothbehaviour}). Our first result is the estimates for some terms of $u(z)$ near the singularity. We improve the estimate of the mixed derivatives when the order of $u$ is $\alpha=1$ and obtain the estimate for higher order derivatives near the origin. We show in \cite{zhang3} that, our result is sharp by use of the generalized hyperbolic metric $\lambda_{\alpha, \beta, \gamma}$ on the thrice-punctured sphere $\mathbb{P}\backslash \{z_1, z_2, z_3\}$ with singularities of order $\alpha, \beta, \gamma \leq 1$ at $z_1, z_2, z_3$, which was given by Kraus, Roth and Sugawa for $\alpha+\beta+\gamma >2$, see \cite{Rothhyper}.
\par As an extremal case of the SK-metric, the hyperbolic metric, also called the Poincar{\'e} metric, plays an important role on (punctured) disks. Early in 1997, Minda \cite{Mindametric} investigated the behavior of the density of the hyperbolic metric in a neighborhood of a puncture on the plane domain using the uniformization theorem. His method offers us a way to describe the asymptotic behavior on an arbitrary hyperbolic region. Our second result is to extend Minda's work and to give some limits of Minda's type.
\par This paper is divided into four sections. In Section 2 the notations and the definitions are introduced. Section 3 is contributed to potential theory. The main results and their proofs are given in Section 4.
\vspace{4mm}

\section {Preliminaries }
\subsection{Singularities and orders}
\setcounter{equation}{0}
\vspace*{3mm}
If $G \subseteq \mathbb{C}$ is a domain, then every positive, upper semi-continuous, real-valued function $\lambda: G \rightarrow (0, + \infty)$ on $G$ induces a conformal metric $\lambda(z)|dz|$, see \cite{Heins,Roth1}. In our discussion we take the linear notation for a conformal metric $ds=\lambda(z)|dz|$. Let $\mathbb{P}$ denote the Riemann sphere $\mathbb{C}\cup \{\infty\}$ and let $\Omega\subseteq\mathbb{P}$ be a subdomain. For a point $p \in \Omega$, let $z$ be local coordinates such that $z(p)=0$. We say a conformal metric $\lambda(z)|dz|$ on the punctured domain $\Omega^*:=\Omega \backslash \{p\}$ has a singularity of order $\alpha\leq 1$ at the point $p$, if, in local coordinates $z$,
\be \label{singularity}
\log\lambda(z)=\left\{
\begin{array}{ll}
-\alpha\log|z|+v(z) & \mbox{if\ }\ \alpha < 1 \\
-\log|z|-\log\log(1/|z|)+w(z)&\mbox{if\ }\ \alpha=1,\end{array} \right.
\ee
where $v(z), w(z) = \textit {O}(1)$ as $z(p)\rightarrow 0$ with $\textit {O}$ and $\textit {o}$ being the Landau symbols throughout our study. Let $M_u(r):=\sup_{|z|=r}u(z)$ for a real-valued function $u(z)$ defined in a punctured neighborhood of $z=0$ and call
\be\label{order of u}
\alpha(u):=\lim_{r\rightarrow 0^+}\frac{M_u(r)}{\log(1/r)}
\ee
the order of $u(z)$ if this limit exists. For $u(z):=\log \lambda(z)$, $\alpha(u)$ in (\ref{singularity}) is equal to $\alpha$ in \eqref{order of u}. In fact, if $\alpha(u)\leq 1$ in \eqref{order of u}, then $v(z)$ is continuous at $z=0$ and $w(z)=\textit {O}(1)$ as $z \rightarrow 0$, see Theorem 3.1 in \cite{Rothbehaviour}. We call the point $p$ a conical singularity or corner of order $\alpha$ if $\alpha< 1$ and a cusp if $\alpha=1$. The generalized Gaussian curvature $\kappa_{\lambda}(z)$ of the density function $\lambda(z)$ is defined by
$$\kappa_{\lambda}(z)=-\frac{1}{\lambda(z)^2}{\liminf_{r\rightarrow 0}\frac{4}{r^2}\left(\frac 1 {2\pi}\int_0^{2\pi}\log\lambda(z+re^{it})dt-\log\lambda(z)\right)}.$$
We say a conformal metric $\lambda(z)|dz|$ on a domain $\Omega\subseteq\mathbb{C}$ is regular, if its density $\lambda(z)$ is positive and twice continuously differentiable, i.e. $\lambda(z)>0$ and $\lambda(z) \in C^2(\Omega)$. If $\lambda(z)|dz|$ is a regular conformal metric, then $$\kappa_{\lambda}(z)=-\frac{\Delta\log\lambda(z)}{\lambda(z)^2},$$
where $\Delta$ denotes the Laplace operator. It is well known that, if $a<\kappa_{\lambda}(z)<b<0$ with constants $a, b \in \mathbb{R}$, the metric $\lambda(z)|dz|$ only has corners or cusps at isolated singularities (see \cite{McOwen1}).
\par The Gaussian curvature is a conformal invariant. Let $\lambda(z)|dz|$ be a conformal metric on a domain $G \in \mathbb{C}$ and $f: \Omega \rightarrow G$ be a holomorphic mapping of a Riemann surface $\Omega$ into $G$. Then we can define the pullback $f^{*}\lambda(w)|dw|$ of $\lambda(z)|dz|$ by
\be
f^{*}\lambda(w)|dw|:=\lambda(f(w))|f'(w)||dw|. \nn
\ee
It is easy to see that $f^{*}\lambda(w)|dw|$ is a conformal metric on $\Omega \backslash\{\mbox{critical points of}\ f\}$ with Gaussian curvature
\be
\kappa_{f^{*}{\lambda}}(w)=\kappa_{\lambda}(f(w)). \nn
\ee
Using this conformal invariance, we can easily build relations between Riemann surfaces with conformal metrics. Here we can see that, on the punctured domain $\Omega \backslash\{\mbox{critical points}$ $\mbox{of}\ f\}$, the critical points of $f$ are the source of the singularities.
\par The hyperbolic metric is a complete metric with some constant Gaussian curvature, here we take the constant to be  $-4$. We call an upper semi-continuous metric $\lambda(z)|dz|$ on a Riemann surface $\Omega$ an SK-metric if its Gaussian curvature is bounded above by $-4$ at every $z\in \Omega$. The hyperbolic metric on the unit disk $\mathbb{D}$ is defined by
\be \label{hyperbolic metric}
\lambda_{\mathbb{D}}(z)|dz|=\frac {|dz|}{1-{|z|}^2}.
\ee
The following result is a fundamental theorem about SK-metrics by Ahlfors, see \cite{Ahlforslemma}, also \cite{Heins}, which claims that the hyperbolic metric $\lambda_{\mathbb{D}}(z)|dz|$ on the unit disk $\mathbb{D}$ is the unique maximal SK-metric on $\mathbb{D}$.

\begin{theorema} $\mathrm{[1]}.$ \label{Ahlfors lemma}
\textsl{Let $ds$ be the hyperbolic metric on $\mathbb{D}$ given in (\ref{hyperbolic metric}) and $d\ell$ be the metrics on $\mathbb{D}$ induce by an SK-metric on a Riemann surface $\Omega$. If the function $f(z)$ is analytic in $\mathbb{D}$, then the inequality
\be
d\ell \leq ds \nn
\ee
will hold throughout the circle. }
\end{theorema}

On the punctured unit disk $\mathbb{D}^*:=\mathbb{D} \backslash \{0\}$, the hyperbolic metric is expressed by
$$\lambda_{\mathbb{D}^*}(z)|dz|=\frac{|dz|}{2|z|\log(1/|z|)}$$
with the constant curvature $-4$. We denote $\mathbb{D}_R:=\{z \in {\mathbb{C}}:|z|<R\}$ and ${\mathbb{D}_R}^*:=\mathbb{D}_R\backslash \{0\}$ for $R>0$. On the punctured disk $\mathbb{D}_R^*$, the (generalized) hyperbolic metric with a conical singularity at the origin is given in \cite{Rothhyper}. For its detailed proof, see \cite{zhang2}.
\begin{theorema} $\mathrm{[7,14]}.$ \label{maximal}
\textsl{ For $R>0$, let
\ben
\lambda_{\alpha,R}(z):=\left\{
\begin{array}{ll}
\displaystyle \frac{(1-\alpha)R^{1-\alpha}|z|^{-\alpha}}{R^{2(1-\alpha)}-|z|^{2(1-\alpha)}} =\frac{1-\alpha}{2|z|\sinh \left((1-\alpha)\log ({R}/{|z|})\right)} & \mbox{if\ }\ \alpha<1, \\
\displaystyle \frac{1}{2|z|\log ({R}/{|z|})} &\mbox{if\ }\ \alpha=1
\end{array} \right.
\een
for $z \in \mathbb{D}_R^*$. Then given an arbitrary SK-metric $\sigma(z)$ on $\mathbb{D}_R^*$ with a singularity at $z=0$ of order $\alpha$, we have $\sigma(z)\leq\lambda_{\alpha,R}(z)$.}
\end{theorema}

\subsection{Regularity and Logarithmic potential}
If equation \eqref{general equation} has a $C^2$-solution $u(z)$, then the higher regularity properties of $u(z)$ only depends on the smoothness of $\kappa(z)$, according to Gilbarg and Trudinger [2, p.\,109]. Here we need the H\"{o}lder spaces $C^{n, \,\nu}(\mathbb{D}_R)$, consisting of functions whose $n$-th order partial derivatives are locally H\"{o}lder continuous with exponent $\nu$ in $\mathbb{D}_R$, $0<\nu \leq 1$, which are defined as the subspaces of $C^n(\mathbb{D}_R)$. The following result can be obtained immediately from the standard regularity theorem, see, e.g. [2,\,Theorem 6.17].
\begin{lemma} $\mathrm{(Regularity\ theorem)}$ \label{regularity coro}
\textsl{ Let $u$ be a $C^2$-solution to the equation $\Delta u=-\kz e^{2u}$ in $\mathbb{D}^*$, where $\kappa \in C^{n,\,\nu}(\mathbb{D}^*)$. Then $u \in C^{n+2,\,\nu}(\mathbb{D}^*)$. If $\kappa$ lies in $C^{\infty}(\mathbb{D}^*)$, then $u \in C^{\infty}(\mathbb{D}^*)$.}
\end{lemma}

\vspace*{2mm}
We shall use potential theory as employed by Kraus and Roth in \cite{Rothbehaviour}. Here we list some elementary facts without proof. \vspace*{2mm}

For a bounded, integrable function $f(z)$ defined on a domain $\Omega \subseteq \mathbb{C}$, the integral
$$\frac {1}{2\pi}\int _{\Omega}L(z-\zeta) f(\zeta) d\sigma_{\zeta}$$
is called the logarithmic potential of $f$, where $L(z-\zeta)=\log|z-\zeta|$ and $d\sigma_{\zeta}$ is the area element on domain $\Omega$. Write $z=x_1+ix_2$, $\zeta=y_1+iy_2$ and set $0<r\leq 1$. The following lemma was mentioned in \cite{Rothbehaviour}. It is a consequence of the famous Riesz decomposition theorem, and can be obtained from Theorem 4.5.1 and Exercise 3.7.3 in \cite{Ransford}.

\begin{lemmaa} $\mathrm{[6]}.$ \label{poisson jensen}
\textsl{Let $u$ be a subharmonic function on $\mathbb{D}_r$ such that $u\in C^2({\mathbb{D}_r}^*)$, $\Delta u$ is integrable in $\mathbb{D}_r$ and
$$\lim_{r\rightarrow 0}\frac{\sup_{|z|=r}u(z)}{\log(1/r)}=0.$$
Then $u(z)=h(z)+\omega(z)$ for $z\in \mathbb{D}_r$, where $h$ is a harmonic function on $\mathbb{D}_r$ and $\omega(z)$ is the logarithmic potential of $\Delta u$.}
\end{lemmaa}

\begin{lemmaa} $\mathrm{[2,\,p.\:\!54]},$ \label{newton potential}
\textsl{ Let $f: \mathbb{D}_r\rightarrow\mathbb{R}$ be a locally bounded, integrable function in $\mathbb{D}_r$ and $\omega$ be the logarithmic potential of $f$. Then $\omega \in C^1({\mathbb{D}_r})$ and for any $z=x_1+ix_2 \in \mathbb{D}_r$,
$$\frac{\pa\:\!{\omega}}{\pa x_j}(z)=\frac 1 {2\pi}\int_{\mathbb{D}_r}\frac {\pa L}{\pa x_j}(z-\zeta) f(\zeta) d\sigma_{\zeta}$$ for $j\ \in \{1,2\}$.\\
If, in addition, $f$ is locally H\"{o}lder continuous with exponent $\nu \leq 1$, then $\omega\in C^2(\mathbb{D}_r)$ and for $z \in \mathbb{D}_r$,
\be
&&\frac{\pa^2\:\!\omega}{\pa x_l\pa x_j}(z)=\frac 1 {2\pi}\int_{\mathbb{D}_R}\frac{\pa^2 L}{\pa x_l\pa x_j}(z-\zeta)\left(f(\zeta)-f(z)\right) d\sigma_{\zeta} \nn\\
&& \qquad \qquad \quad \: \ \,-\frac{1}{2\pi}f(z)\int_{\pa \mathbb{D}_R}\frac{\pa L}{\pa x_j}(z-\zeta)N_l(\zeta)|d\zeta|, \nn
\ee
where $N(\zeta)=(N_1(\zeta),N_2(\zeta))$ is the unit outward normal at the point $\zeta \in\pa \mathbb{D}_R$ with $R>r$ and $f$ is extended to vanish outside of $\mathbb{D}_r$.}
\end{lemmaa}
\par There is a similar proposition for higher order derivatives of the logarithmic potential. Define a multi-index $\mj=(j_1, j_2)$,  $|\mj|=j_1+j_2$, $j_1, j_2=0,1,2,\ldots \,$, so $(\zeta-z)^{\mj}=(y_1-x_1)^{j_1}(y_2-x_2)^{j_2}$, $\mj !=j_1!j_2!$. For $z=x_1+ix_2$, denote $$\frac{\pa}{\pa x_1}=\pa_1,\ \frac{\pa}{\pa x_2}=\pa_2,\  \displaystyle \pa^{\mj}=\pa^{j_1}_1 \pa^{j_2}_2.$$ For a given multi-index $\mj=(j_1, j_2)$, we can choose $\me_{\tau}=(0,1)$ or $(1,0)$ for $\tau=1,2,\ldots \,$ such that $\mj=\me_1+\me_2+\cdots+\me_n$ with $n=|\mj|$. Write $\zeta=y_1+iy_2$, set \ben
P_n[f](z,\zeta):=\left\{
\begin{array}{ll}
\vspace*{2mm}
\displaystyle \sum_{|\ma|\leq n} \frac{(\zeta-z)^{\ma}{\pa}^{\ma}}{\ma !}f(z) & \mbox{if\ }\ n \geq 1 \\
f(z) & \mbox{if\ }\ n=0,
\end{array} \right.
\een
where $\ma$ is a multi-index. For $m=1,2$, it holds that $$\frac{\pa P_n[f]}{\pa y_m}(z,\zeta)=P_{n-1}[\pa_m f](z,\zeta),$$
see \cite{zhang2} for more details. Using this notation, we can present the analogue of Lemma \ref{newton potential} as follows.
\begin{lemmaa} $\mathrm{[14]}.$ \label{new higher prop}
\textsl{ Let $0<r<1$, $f: \mathbb{D}_r\rightarrow\mathbb{R}$ and $f \in C^{n-2,\,\nu}(\mathbb{D}_r)$ with $0<\nu \leq 1$, $n\geq 3$, $\omega$ be the logarithmic potential of $f$. Then
$\omega(z)\in C^{n}(\mathbb{D}_r)$ and for $n=|\mj|$,
\be \label{new higher}
&& \pa^{\mj}\omega(z) \nn \\
&=& \frac 1 {2\pi}\int_{\mathbb{D}_R}\pa^{\mj} L(z-\zeta) \cdot \left(f(\zeta)-P_{n-2}[f](z,\zeta)\right)d\sigma_{\zeta} \nn \\
&& -\frac{1}{2\pi} \sum^{n-1}_{\tau=1} \int_{\pa \mathbb{D}_R}\pa^{\mtheta_{\tau}} L(z-\zeta) \cdot P_{\tau-1}[\pa^{\mphi_{\tau}}f](z,\zeta)\cdot \langle N(\zeta),\me_{\tau+1} \rangle |d\zeta|,
\ee
where $\mtheta_{\tau}:=\me_1+ \cdots +\me_{\tau}$, $\mphi_{\tau}:=\me_{\tau+2}+ \cdots +\me_n$ for $\tau=1, \ldots, \,n-1$ and $\mphi_{n-1}:=(0,0)$, $N(\zeta)=(N_1(\zeta),N_2(\zeta))$ is the unit outward normal at the point $\zeta \in \pa \mathbb{D}_R$ with $R>r$, $\langle \ , \  \rangle$ is the inner product and the function $f$ is extended to vanish outside of $\mathbb{D}_r$.}
\end{lemmaa}

\section{Main estimates}
\setcounter{equation}{0}
We denote $$\pa^n=\frac{\pa ^n}{\pa z^n},\ \bar{\pa}^n=\frac {\pa^n}{\pa\bar{z}^n}$$ for $n\geq1$. The following theorem is given by Kraus and Roth in \cite{Rothbehaviour}.
\begin{theorema} $\mathrm{[6]}.$ \label{original}
\textsl{ Let $\kappa:\mathbb{D}\rightarrow\mathbb{R}$ be a locally H\"{o}lder continuous function with $\kappa(0)<0$. If $u:\mathbb{D}^*\rightarrow \mathbb{R}$ is a $C^2$-solution to $\Delta u=-\kz e^{2u}$ in $\mathbb{D}^*$, then $u$ has an order $\alpha \in (-\infty,1]$ and
\begin{align}
&u(z)=-\alpha\log|z|+v(z),            & \textrm{if\ } \ \alpha<1,\nonumber \\
&u(z)=-\log|z|-\log\log(1/|z|)+w(z),  & \textrm{if\ } \ \alpha=1,\nonumber
\end{align}
where the remainder functions $v(z)$ and $w(z)$ are continuous in $\mathbb{D}$. Moreover, the first partial derivatives with respect to $z$ and $\bar{z}$,
\begin{align}
& \pa v(z),\,\bar{\pa} v(z)\ \mbox{are continuous at}\ z=0& \mbox{if\ }&\ \alpha<1/2;\nonumber
\end{align}
and
\begin{align}
& \pa v(z),\,\bar{\pa} v(z)=\textit {O}(1) &\mbox{if\ }& \ \alpha=1/2;\nonumber\\
& \pa v(z),\,\bar{\pa} v(z)=\textit {O}(|z|^{1-2\alpha}) &\mbox{if\ }&\  1/2<\alpha<1,\nonumber\\
& \pa w(z),\,\bar{\pa} w(z)=\textit {O}(|z|^{-1}(\log(1/|z|))^{-2}) &\mbox{if\ }& \ \alpha=1,\nonumber
\end{align}
when $z$ approaches $0$. In addition, the second partial derivatives,
\begin{align}
&\pa^2 v(z),\,\pa \bar{\pa}v(z) \ \textrm{and}\  \bar{\pa}^2 v(z) \ \textrm{are\  continuous \ at}\ z=0 & \textrm{if\ }& \ \alpha \leq 0;\nonumber\end{align}
and
\begin{align}
& \pa^2 v(z),\,\pa \bar{\pa}v(z),\,\bar{\pa}^2 v(z)= \textit {O}(|z|^{-2\alpha}) &\textrm{if\ }&\ 0<\alpha<1,\nn\\
& \pa^2 w(z),\,\pa \bar{\pa} w(z),\,\bar{\pa}^2 w(z)=\textit {O}(|z|^{-2}(\log(1/|z|))^{-2}) & \textrm{if\ }& \ \alpha=1,
\end{align}
when $z$ tends to $z=0$.}
\end{theorema}
\par In the work of Kraus and Roth, the proof of Theorem \ref{original} was based on Lemma \ref{newton potential}. Since we have obtained a similar statement in Lemma \ref{new higher prop}, the estimate for higher order derivatives of the remainder functions $v(z)$, $w(z)$ can be given. We consider $v(z)$ and $w(z)$ separately.
\begin{theorem} \label{estimate v}
\textsl{ Let $\kappa(z)$, $u(z)$, $v(z)$ and $\alpha$ be the same as in Theorem \ref{original}. If $0<\alpha<1$ and if, in addition, $\kz \in C^{n-2,\,\nu}(\mathbb{D}^*)$ for an integer $n \geq 3$, $0<\nu \leq 1$, then for $n_1$, $n_2 \geq 1$, $n_1+n_2 =n$, near the origin, the remainder function $v(z)$ satisfies
$$\pa^n v(z),\,\bar{\pa}^n  v(z),\,\bar{\pa}^{n_1}\pa^{n_2}v(z)=\textit{O}(|z|^{2-2\alpha-n}).$$}
\end{theorem}
\textbf{Proof.}  Lemma \ref{regularity coro} shows that $u(z) \in C^{n,\,\nu}(\mathbb{D}^*)$. Due to Kraus and Roth \cite{Rothbehaviour} we have $$v(z)=h(z)+\frac 1 {2\pi}\int_{\mathbb{D}_r}L(z-\zeta)f(\zeta)
d\sigma_{\zeta}$$ for $z\in \mathbb{D}_r^*$, $0<r<1$ and a harmonic function $h$ on $\mathbb{D}_r$, where $q(z)=-\kappa(z)e^{2v(z)}$, $f(z)=q(z)|z|^{-2\alpha}$. Now fix $0<R<1$, choose $z \in \mathbb{D}_{R/2}^*$ and let $r=|z|/2$. Then for a multi-index $\mj$, $|\mj|=n \geq 3$, rearranging (\ref{new higher}) leads to
\be \label{for v}
&&\pa^{\mj}v(z) \nn \\
&=&\pa^{\mj}h(z)+\frac{1}{2\pi}\int_{\mathbb{D}_{R}\backslash \mathbb{D}_{r}}\pa^{\mj}L(z-\zeta)f(\zeta)d\sigma_{\zeta}+  \frac 1 {2\pi}\int_{\mathbb{D}_r}\pa^{\mj}L(z-\zeta)
\left( f(\zeta)-f(z)\right)d\sigma_{\zeta} \nonumber\\
&&+\frac 1 {2\pi}\int_{\mathbb{D}_r}\pa^{\mj}L(z-\zeta)\sum_{1\leq |\ma|\leq n} \frac{(\zeta-z)^{\ma}\pa^{\ma}f(z)}{\ma!}d\sigma_{\zeta} \nonumber\\
&&-\frac{1}{2\pi} \sum^{n-1}_{\tau=1} \int_{\pa \mathbb{D}_r}\pa^{\mtheta_{\tau}} L(z-\zeta) \cdot P_{\tau-1}[\pa^{\mphi_{\tau}}f](z,\zeta)\cdot \langle N(\zeta),\me_{\tau+1} \rangle |d\zeta|
\ee
for $z=x_1+ix_2$ and a harmonic function $h$ on $\mathbb{D}_R$, and the same symbols $\mtheta_{\tau}$, $\mphi_{\tau}$ are used here as in \eqref{new higher}.

It is known that
\be \label{logn}
\left| \pa^{\mj}L(z-\zeta) \right|
\leq \frac{n!}{|z-\zeta|^{n}},
\ee
see [2, p.\,17]. Denote $M=\sup_{\zeta \in \mathbb{D}_R}|q(\zeta)|$ and let $C_n>0$, $n \in \mathbb{N}$, be some constants. Then
$$\left| \int_{\mathbb{D}_{R}\backslash \mathbb{D}_{r}}\pa^{\mj}L(z-\zeta)f(\zeta)d\sigma_{\zeta} \right| \leq  M \int_{\mathbb{D}_{R}\backslash \mathbb{D}_{r}}
\frac{n!}{|z-\zeta|^{n}} \frac{1}{|\zeta|^{2\alpha}} d\sigma_{\zeta} \leq
\frac{C_1}{|z|^{2\alpha+n-2}},$$
and
\be
&&\left| \int_{\mathbb{D}_r}\pa^{\mj}L(z-\zeta)
\left( f(\zeta)-f(z)\right)d\sigma_{\zeta} \right| \nn \\
&\leq&  \int_{\mathbb{D}_{r}} \frac{n!}{|z-\zeta|^{n}} \frac{|q(\zeta)-q(z)|}{|\zeta|^{2\alpha}} d\sigma_{\zeta}
+ M \int_{\mathbb{D}_{r}}\frac{n!}{|z-\zeta|^{n}}
\frac{(|\zeta|^{\alpha}+|z|^{\alpha})||\zeta|^{\alpha}-|z|^{\alpha}|}{|z|^{2\alpha}|\zeta|^{2\alpha}} d\sigma_{\zeta} \nn \\
&\leq& \frac{C_2}{|z|^{2\alpha+n-2}}, \nn
\ee
see \cite{Rothbehaviour}. When one of $j_1$ and $j_2$ is zero, there is no cancelation. Thus we obtain
$\pa^n v,  \bar{\pa}^n  v=\textit{O}(|z|^{2-2\alpha-n}).$ If neither $j_1$ nor $j_2$ is zero, the first three integrals in (\ref{for v}) are canceled, so we have to consider the last term in (\ref{for v}). In the last sum in (\ref{for v}) letting $\tau=1$, we get the integral
$$\pa^{\mphi_1}f(z) \cdot \int_{\pa \mathbb{D}_r}\pa^{\me_1} L(z-\zeta) \cdot \langle N(\zeta),\me_{2} \rangle |d\zeta|. $$
Writing $\zeta=re^{i\theta}$ and taking $\me_1=(1,0)$, $\me_2=(0,1)$ without loss of generality, we have
\be \label{sincos}
&&\left|\int_{\pa \mathbb{D}_r}\pa_1L(z-\zeta)N_2(\zeta)|d\zeta|\right|=\left|\int_{\pa \mathbb{D}_r}\frac{x_1-r\cos\theta}{|z-\zeta|^2}\sin\theta|d\zeta|\right|\nn\\
&=&\left|\int^{2\pi}_0 \frac{x_1-r\cos\theta}{|z-\zeta|^2} r \sin\theta  d\theta\right| \leq 2\pi \frac{|x_1-r\cos\theta|}{|z-\zeta|^2} r |\sin\theta| \leq 6\pi.
\ee
So it is evident that,
$$\left|\int_{\pa \mathbb{D}_r}\pa^{\me_1} L(z-\zeta) \cdot \langle N(\zeta),\me_{2} \rangle |d\zeta| \right| \leq 6\pi$$
holds for all kinds of $\me_1$ and $\me_2$. Now consider $\pa^{\mphi_1}f(z)$. Since $f(z)=q(z)|z|^{-2 \alpha}$, then the term $q(z)\cdot \pa^{\mphi_1}(|z|^{-2 \alpha})$ appears with some coefficient. Note that
$$\left|\pa^{\mphi_1} \left(\frac{1}{|z|^{2\alpha}}\right)\right| \leq \frac{C_3}{|z|^{2\alpha+n-2}},$$
so
$$\left| q(z) \pa^{\mphi_1} \left(\frac{1}{|z|^{2\alpha}}\right) \cdot \int_{\pa \mathbb{D}_r}\pa^{\me_1} L(z-\zeta) \cdot \langle N(\zeta),\me_{2} \rangle |d\zeta| \right| \leq \frac{6\pi M C_3}{|z|^{2\alpha+n-2}}.$$ Therefore $\bar{\pa}^{n_1}\pa^{n_2}v=\textit{O}(|z|^{2-2\alpha-n}). $  \hspace*{\fill} $\Box$\\
\par The following result is for the higher order derivatives of the remainder functions $w(z)$ when the order is $1$.
\begin{theorem} \label{for w theorem}
\textsl{ Let $\kappa(z)$, $u(z)$, $w(z)$ and $\alpha$ be the same as in Theorem \ref{original}. If $\alpha=1$ and if, in addition, $\kz \in C^{n-2,\,\nu}(\mathbb{D}^*)$ for an integer $n \geq 3$, $0<\nu \leq 1$, then for $n_1$, $n_2 \geq 1$, $n_1+n_2 =n$, near the origin, the remainder functions $w(z)$ satisfies
\be
\bar{\pa}^nw(z),\,\pa^nw(z)=\textit {O}(|z|^{-n}(\log(1/|z|))^{-2}), \nn \\
\bar{\pa}^{n_1}\pa^{n_2}w(z)=\textit {O}(|z|^{-n}(\log(1/|z|))^{-3}).  \label{new mixed}
\ee}
\end{theorem}
\par The proof is based on the following lemma.
\begin{lemmaa}\label{estimate of kzew-1} $\mathrm{[6]}.$ \ \textsl{Let $\kappa:\mathbb{D}\rightarrow\mathbb{R}$ be a continuous function with $\kappa(0)<0$ and
\be
\kz=\kappa(0)+\frac{s(z)}{(\log(1/|z|))^2},\nonumber\ee
where $s(z)=\textit {O}(1)$ as $z\rightarrow 0$. If $u: \mathbb{D}^* \rightarrow\mathbb{R}$ is a solution to $\Delta u=-\kz e^{2u}$ with $u(z)=-\log|z|-\log\log(1/|z|)+w(z)$ where $w(z)=\textit {O}(1)$ for $z\rightarrow0$, then there exists a disk $\mathbb{D}_{\rho}$ such that
\be \label{lemma}
\left|-\kz e^{2w(z)}-1\right|\leq\frac C {\log(1/|z|)},\ \ \ z\in \mathbb{D}_{\rho},
\ee
for some constant $C>0$.}
\end{lemmaa}
\textbf{Proof of Theorem \ref{for w theorem}.} Lemma \ref{regularity coro} shows that $u(z) \in C^{n,\,\nu}(\mathbb{D}^*)$.
For $w(z)$ defined as in Theorem \ref{original}, we first show that
\be\label{w(z)}
w(z)=h(z)+\frac 1 {2\pi}\int_{\mathbb{D}_r}L(z-\zeta)\frac{-\kappa(\zeta)e^{2w(\zeta)}-1}{|\zeta|^2(\log(1/|\zeta|))^2}
d\sigma_{\zeta}\ee
for $z\in \mathbb{D}_r^*$, $0<r<1$, where $h$ is harmonic on $\mathbb{D}_r$. Let
$$t(z):=-\log\log(1/|z|), \quad p(z):=w(z)+t(z)=u(z)+\log|z|$$
for $z\in \mathbb{D}_r^*$. Since
$$\Delta p(z)=-\kz e^{2u}=\frac{-\kz e^{2w(z)}}{|z|^2(\log(1/|z|))^2}>0,$$
$p(z)$ is subharmonic on $\mathbb{D}_r^*$ and $\lim_{z\rightarrow0}p(z)=-\infty$, then $p(z)$ is subharmonic on $\mathbb{D}_r$. By Lemma \ref{poisson jensen}, as $z \mapsto \Delta p(z)$ is integrable over $\mathbb{D}_r$,
\be
p(z)=h_p(z)+\frac 1 {2\pi}\int_{\mathbb{D}_r}L(z-\zeta)\frac{-\kappa(\zeta)e^{2w(\zeta)}}{|\zeta|^2(\log(1/|\zeta|))^2}
d\sigma_{\zeta},\ z\in \mathbb{D}_r,\nonumber\ee where $h_p(z)$ is harmonic on $\mathbb{D}_r$.
For $t(z)$, we also have
\be
t(z)=h_t(z)+\frac 1 {2\pi}\int_{\mathbb{D}_r}L(z-\zeta)\frac{1}{|\zeta|^2(\log(1/|\zeta|))^2}d\sigma_{\zeta},\ z\in \mathbb{D}_r,\nonumber\ee where $h_t(z)$ is harmonic on $\mathbb{D}_r$. Setting $w(z)=p(z)-t(z)$ gives \eqref{w(z)} with $h(z)=h_p(z)-h_t(z)$.

Now set $R<1 / e^2$. So there exists a number ${\rho}>0$ such that the inequality (\ref{lemma}) holds in the disk $\mathbb{D}_{\rho}$. Let $\widetilde{\rho}=\min\{R/2,\,\rho\}$. We choose $z\in \mathbb{D}_{\widetilde{\rho}}$ and set $r=|z|/2$. Let $q(z)=-\kappa(z)e^{2w(z)}-1$, $f(z)=q(z)|z|^{-2\alpha}$. Then from (\ref{new higher}), we have
\be \label{for w}
&& \pa^{\mj}w(z)\nn \\
&=&\pa^{\mj}h(z)+\frac{1}{2\pi}\int_{\mathbb{D}_{\widetilde{\rho}}\backslash \mathbb{D}_{r}}\pa^{\mj}L(z-\zeta)f(\zeta)d\sigma_{\zeta}+  \frac 1 {2\pi}\int_{\mathbb{D}_r}\pa^{\mj}L(z-\zeta)
\left( f(\zeta)-f(z)\right)d\sigma_{\zeta} \nn\\
&&+\frac 1 {2\pi}\int_{\mathbb{D}_r}\pa^{\mj}L(z-\zeta)\sum_{1\leq |\ma|\leq n} \frac{(\zeta-z)^{\ma}\pa^{\ma}f(z)}{\ma!}\;d\sigma_{\zeta} \nonumber\\
&&-\frac{1}{2\pi} \sum^{n-1}_{\tau=1} \int_{\pa \mathbb{D}_r}\pa^{\mtheta_{\tau}} L(z-\zeta) \cdot P_{\tau-1}[\pa^{\mphi_{\tau}}f](z,\zeta)\cdot \langle N(\zeta),\me_{\tau+1} \rangle |d\zeta|
\ee
for a harmonic function $h$ on $\mathbb{D}_{\widetilde{\rho}}$. We can obtain
$$\left| \int_{\mathbb{D}_{\widetilde{\rho}}\backslash \mathbb{D}_{r}}\pa^{\mj}L(z-\zeta)f(\zeta)d\sigma_{\zeta} \right| \leq
\frac{C_4}{|z|^n (\log(1/|z|))^2},$$
$$ \left| \int_{\mathbb{D}_r}\pa^{\mj}L(z-\zeta)
\left( f(\zeta)-f(z)\right)d\sigma_{\zeta} \right| \leq \frac{C_5}{|z|^n (\log(1/|z|))^2},$$
by \eqref{logn} and Theorem \ref{logn} in \cite{Rothbehaviour}. So $ \bar{\pa}^nw(z)$, $\pa^nw(z)=\textit {O}(|z|^{-n}(\log(1/|z|))^{-2})$. For the mixed partial derivatives, since the first three integrals are canceled, we have to estimate the last term in (\ref{for w}). Letting $\tau=1$ in the last sum of (\ref{for w}), the term
$$ \pa^{\mphi_{1}}f(z)\int_{\pa \mathbb{D}_r}\pa^{\me_1} L(z-\zeta) \cdot \langle N(\zeta),\me_{2} \rangle |d\zeta|$$
appears. Now consider $\pa^{\mphi_{1}}f(z)$. Our aim is $\bar{\pa}^{n_1}\pa^{n_2}w(z)=\textit {O}(|z|^{-n}(\log(1/|z|))^{-3})$. Since $f(z)=q(z)|z|^{-2}(\log(1/|z|))^{-2}$, then $q(z)\cdot \pa^{\mphi_{1}} (|z|^{-2}(\log(1/|z|))^{-2})$ appears in $ \pa^{\mphi_{1}} f(z)$ with some coefficient. We can calculate that
$$\left| \pa^{\mphi_{1}} \frac{1}{|z|^{2}(\log(1/|z|))^{2}}\right| \leq \frac{C_6}{|z|^{n}(\log(1/|z|))^2},$$
thus
$$\left|  q(z) \pa^{\mphi_{1}} \frac{1}{|z|^{2}(\log(1/|z|))^{2}} \cdot \pa^{\me_1} L(z-\zeta) \cdot \langle N(\zeta),\me_{2} \rangle |d\zeta|\right| \leq \frac{6\pi C_7 \cdot C_6}{|z|^{n}(\log(1/|z|))^{3}}$$
provided (\ref{logn}) and (\ref{sincos}). So $ \bar{\pa}^{n_1}\pa^{n_2}w(z)=\textit {O}(|z|^{-n}(\log(1/|z|))^{-3})$ as desired. \hspace{\fill} $\Box$\\

The second order derivative of $w(z)$ in Theorem \ref{original} is contained in Theorem \ref{for w theorem}. However, for the mixed partial derivative, the estimate (\ref{new mixed}) is more accurate than (4.1). We take it to be a corollary as following.
\begin{corollary} \label{second partial derivative order}
\textsl{ Let $\kappa:\mathbb{D}\rightarrow\mathbb{R}$ be a locally H\"{o}lder continuous function with $\kappa(0)<0$. If $u:\mathbb{D}^* \rightarrow \mathbb{R}$ is a $C^2$-solution to $\Delta u=-\kz e^{2u}$ in $\mathbb{D}^*$ with the order $\alpha=1$ at the point $z=0$, then for the remainder function $w(z)$ as in Theorem \ref{original}, the second partial derivatives satisfy
$$w_{z\bar{z}}(z)=\textit {O}(|z|^{-2}(\log(1/|z|))^{-3}).$$ }
\end{corollary}

As for the sharpness of Theorems \ref{original}, \ref{estimate v} and \ref{for w theorem}, the generalized hyperbolic metric on the thrice-punctured sphere makes a convictive case here. Theorems 3.3 and 4.2 in \cite{zhang3} verify that Theorems \ref{original}, \ref{estimate v} and \ref{for w theorem} are sharp, see \cite{zhang3} for details.

\vspace*{3mm}

\section{Minda-type theorems}
\setcounter{equation}{0}
The following result is Minda's theorem. It is a general estimate for the hyperbolic metric near the singularity.
\begin{theorema} $\mathrm{[10]}.$ \label{Minda original}
\textsl{Suppose $\Omega$ is a hyperbolic region in the complex plane and $p \in \mathbb{C}$ is an isolated boundary point of $\Omega$. Let the hyperbolic metric on $\Omega$ with the constant Gaussian curvature $-1$ be $\lambda_{\Omega}(\omega)|d\omega|$. Then
\ben
\lim_{\omega\rightarrow p}|\omega-p|\log(1/|\omega-p|)\lambda_{\Omega}(\omega)=\frac{1}{2}\,.
\een}
\end{theorema}

The following theorem is due to Kraus and Roth.
\begin{theorema} $\mathrm{[6]}.$ \label{estimate for cusps original}
\textsl{Let $\lambda(z)|dz|$ be a regular conformal metric on a domain $\Omega\subseteq\mathbb{C}$ with an isolated singularity at $z=p$. Suppose that its curvature $\kappa :\Omega\rightarrow \mathbb{R}$ has a H\"{o}lder continuous extension to $\Omega\cup\{p\}$ such that $\kappa(p)< 0$. Then $\log\lambda$ has an order $\alpha\leq1$ at $z=p$ and
\ben
\lim_{z\rightarrow p}|z-p|\log(1/|z-p|)\lambda(z)=\left\{\begin{array}{ll}
0&\mbox{if\ }\ \alpha< 1\\
\displaystyle \frac 1 {\sqrt{-\kappa(p)}}&\mbox{if\ }\ \alpha=1.\end{array}\right.
\een
\vspace*{1mm} }
\end{theorema}
\par We obtain the following result in relation to Theorem \ref{estimate for cusps original}.
\begin{theorem}\label{coro}
\textsl{ Let $\lambda(z)|dz|$ be a regular conformal metric on a domain $\Omega\subseteq\mathbb{C}$ with an isolated singularity at $z=p$. Suppose that the curvature $\kappa :\Omega\rightarrow \mathbb{R}$ has a H\"{o}lder continuous extension to $\Omega\cup\{p\}$ such that $\kappa(p)< 0$ and the order of $\log\lambda$ is $\alpha=1$ at $z=p$. Then \vspace*{2mm} \\
(i)$\ \displaystyle \lim_{z\rightarrow p} (z-p)|z-p|\log(1/|z-p|)\lambda_{z}(z)=-{\frac 1 {2\sqrt{-\kappa(p)}}}, \hspace*{\fill} $ \vspace*{2mm}\\
(ii)$\ \displaystyle \ \lim_{z\rightarrow p}(z-p)^2|z-p|\log(1/|z-p|)\lambda_{zz}(z)=\displaystyle {\frac 3{4\sqrt{-\kappa(p)}}}, \hspace*{\fill} $ \vspace*{2mm}\\
(iii)$\ \displaystyle \ \lim_{z\rightarrow p}|z-p|^3\log(1/|z-p|)\lambda_{z\bar{z}}(z)=\displaystyle {\frac 1 {4\sqrt{-\kappa(p)}}}. \hspace*{\fill} $ }
\end{theorem}
\textbf{Proof.} Let $\mathbb{H}$ be the upper half-plane. For each simply closed curve  $\gamma: [0,1] \rightarrow \Omega$ around $p$ with $\gamma(0)=\gamma(1)$, consider the lift $\widetilde{\gamma}$ of $\gamma$ in $\mathbb{H}$. Since there exists an automorphism $g$ on $\mathbb{H}$ such that $\widetilde{\gamma}(1)=g(\widetilde{\gamma}(0))$, we may assume that $g(z)=z+1$ on $\mathbb{H}$. Let $\pi: \mathbb{H}\rightarrow \Omega$ be the regular covering projection, we have $\pi \circ g = \pi$. Define $\varphi: \mathbb{H} \rightarrow \mathbb{D}^*$, $z \mapsto e^{2 \pi iz}$, then the quotient space $\mathbb{H}/ \langle g \rangle$ is conformally equivalent to $\mathbb{D}^*$. Hence there exists a conformal mapping $\rho: \mathbb{D}^* \rightarrow \Omega$ such that $\rho \circ \varphi=\pi$ and $\rho$ can be extended to $\mathbb{D} \rightarrow \Omega$ holomorphically by setting $\rho (0)=p$.
So it is sufficient to consider the case $p=0$ and $\mathbb{D}^*=\Omega$.
\par Let $u(z):=\log\lambda(z)$, so $\lambda_z(z)=u_z(z)\lambda(z)$. It holds
$$\lim_{z\rightarrow 0}z u_z(z)=-\frac 1 2,\; \;
\lim_{z\rightarrow 0}z^2 u_{zz}(z)=\frac 1 2,\; \; \lim_{z\rightarrow0}|z|^2u_{z\bar{z}}=0$$
by Theorem \ref{original}. In combination with Theorem \ref{estimate for cusps original}, we have
\ben
&&\lim_{z\rightarrow 0}z|z|\log(1/|z|)\lambda_z(z)=\lim_{z\rightarrow 0}z|z|\log(1/|z|) u_z(z)\lambda(z)\\
&=&\lim_{z\rightarrow 0}|z|\log(1/|z|)\lambda(z) \cdot z u_{z}(z)=-\frac 1 {2\sqrt{-\kappa(0)}}
\een
for the first case,
\ben
&&\lim_{z\rightarrow 0}z^2|z|\log(1/|z|)\lambda_{zz}(z)=\lim_{z\rightarrow 0}z^2|z|\log(1/|z|)(u_{zz} \lambda+u_z \lambda_z)\\
&=&\lim_{z\rightarrow 0}(z^2u_{zz}\cdot|z|\log(1/|z|)\lambda)+\lim_{z\rightarrow 0}(z|z|\log(1/|z|)\lambda_z\cdot zu_z)\\
&=&\frac 1 {2\sqrt{-\kappa(0)}}+(-\frac 1 2)\cdot(-\frac 1 {2\sqrt{-\kappa(0)}})=\frac 3{4\sqrt{-\kappa(0)}}\nonumber
\een
for the second case and
\ben
&&\lim_{z\rightarrow0}|z|^3\log(1/|z|)\lambda_{z\bar{z}}(z)=\lim_{z\rightarrow0}|z|^3\log(1/|z|)(u_{z\bar{z}}\lambda+u_{z}
\lambda_{\bar{z}})\nonumber\\
&=&\lim_{z\rightarrow0}(|z|^2u_{z\bar{z}}\cdot|z|\log(1/|z|)\lambda)
+\lim_{z\rightarrow0}(\bar{z}|z|\log(1/|z|)\lambda_{\bar{z}}\cdot zu_z)\nn\\
&=&-\frac 1 {2\sqrt{-\kappa(0)}}\cdot(-\frac 1 2)=\frac 1 {4\sqrt{-\kappa(0)}}\nonumber
\een
for the last case as desired. \hfill  $\Box$
\vspace*{2mm}
\par Theorem \ref{coro} is given for a regular conformal metric with a (locally) H\"oder continuous Gaussian curvature $\kappa$. Considering Theorems \ref{for v} and \ref{for w theorem}, if we add the assumption that $\kappa$ is $n$-th order (locally) H\"older continuous, we can obtain the higher order version of Theorems \ref{estimate for cusps original} and \ref{coro}.
\begin{theorem}\label{general u}
\textsl{Let $\kappa:\mathbb{D}\rightarrow\mathbb{R}$ be of class $C^{n-2,\,\nu}(\mathbb{D}^*)$ for an integer $n \geq 3$, $0<\nu \leq 1$ and $\kappa(0)<0$. If $u:\mathbb{D}^*\rightarrow \mathbb{R}$ is a $C^{n,\,\nu}$-solution to $\Delta u=-\kz e^{2u}$ in $\mathbb{D}^*$, then $u$ has order $\alpha \in (-\infty, 1]$ and for $n_1,\,n_2 \geq 1$, $n_1+n_2\leq n$, \vspace*{2mm} \\
(i) $ \displaystyle \ \lim_{z\rightarrow0}z^n\pa^nu(z)=\frac {\alpha}2(-1)^n(n-1)!=\lim_{z\rightarrow0}\bar{z}^n\bar{\pa}^nu(z),$ \vspace*{2mm} \\
(ii) $ \displaystyle \ \lim_{z\rightarrow 0}\bar{z}^{n_1} z^{n_2}\bar{\pa}^{n_1}\pa^{n_2}u(z)=0.$}
\end{theorem}
\textbf{Proof.} When $0<\alpha<1$,\ $u(z)=-\alpha\log|z|+v(z)$. Theorems 4.1 and 4.2 imply that
$$\lim_{z \rightarrow 0}z^n \pa^n v(z)=0, \ \ \lim_{z \rightarrow 0}\bar{z}^{n_1} z^{n_2} \bar{\pa}^{n_1}\pa^{n_2} v(z)=0$$
for $n_1$, $n_2$, $n \geq 1$. Since
\be \label{pa logz}
\pa^{n}\log|z|=\frac{(-1)^{n-1}(n-1)!}{2z^n}, \quad \bar{\pa}^{n_1}\pa^{n_2}\log|z|=0,
\ee
so $$\lim_{z \rightarrow 0}z^n \pa^n u(z)=-\alpha \lim_{z \rightarrow 0}z^n \pa^n \log|z|+\lim_{z \rightarrow 0}z^n \pa^n v(z)=\frac{\alpha}{2z^n}(-1)^{n}(n-1)!,$$
$$\lim_{z \rightarrow 0}\bar{z}^{n_1} z^{n_2} \bar{\pa}^{n_1}\pa^{n_2} u(z)=0.$$
When $\alpha=1$,\ $u(z)=-\log|z|-\log\log(1/|z|)+w(z)$. We have
$$\lim_{z \rightarrow 0}z^n \pa^n w(z)=0, \ \ \lim_{z \rightarrow 0}\bar{z}^{n_1} z^{n_2} \bar{\pa}^{n_1}\pa^{n_2} w(z)=0$$
for $n_1,\,n_2,\,n \geq 1$, from Theorems 4.1 and 4.3. By induction,
$$\pa^n\log\log(1/|z|)=\sum^n_{j=1}\frac{C^{(n)}_j}{z^n(\log(1/|z|))^j}$$
with constant $C^{(n)}_{j}$ for $1 \leq j \leq n$. If we fix $n_2$, then
$$\bar{\pa}^{n_1}\pa^{n_2}\log\log(1/|z|)=\sum^{n_2}_{j=1}\frac{C^{(n_1,\, n_2)}_{j}}{\bar{z}^{n_1}z^{n_2} (\log(1/|z|))^{j+1}}$$
with constant $C^{(n_1,\,n_2)}_{j}$ for $1 \leq j \leq n_2$.
So $$\lim_{z \rightarrow 0}z^n \pa^n \log\log(1/|z|)=0, \ \ \lim_{z \rightarrow 0}\bar{z}^{n_1}z^{n_2} \bar{\pa}^{n_1}\pa^{n_2} \log\log(1/|z|)=0$$ for $n_1,\,n_2,\,n \geq 1$. Combining with \eqref{pa logz} leads to
$$\lim_{z \rightarrow 0}z^n \pa^n u(z)=-\alpha \lim_{z \rightarrow 0}z^n \pa^n \log|z|+\lim_{z \rightarrow 0}z^n \pa^n v(z)=\frac{(-1)^{n}(n-1)!}{2z^n},$$
$$\hspace*{51mm}\lim_{z \rightarrow 0}\bar{z}^{n_1} z^{n_2} \bar{\pa}^{n_1}\pa^{n_2} u(z)=0.\hspace*{51mm} \Box $$

From the proof above, we can obtain a stronger limit for the mixed derivative of $u(z)$ when the order $\alpha=1$,
$$\displaystyle \ \lim_{z\rightarrow 0}\bar{z}^{n_1} z^{n_2}(\log(1/|z|))^2\bar{\pa}^{n_1}\pa^{n_2}u(z)=C^{(n_1,\;n_2)}_1=\frac{(-1)^{n_1+n_2-1}}{4}(n_1-1)!(n_2-1)!,$$
see \cite{zhang3} for more details.
\vspace*{2mm}
\par On the basis of Theorem \ref{general u}, we can provide the following result as a higher order estimate for a conformal metric with the negative curvature near the origin when $\alpha=1$.
\begin{theorem} \label{higher lambda for cusps}
\textsl{Let $\kappa$ and $u$ be the same as in Theorem \ref{general u}. If the order of $u$ is $\alpha=1$, then for $n_1,\,n_2 \geq 0$, $n_1+n_2 \leq n$, the limit
$$l_{n_1,n_2}:=\frac{1}{n_1!n_2!}\lim_{z \rightarrow0}|z|\log(1/|z|){\bar{z}}^{n_1} z^{n_2} {\bar{\pa}}^{n_1} \pa^{n_2} \lambda(z)$$
exists. Moreover, the numbers $l_{n_1,n_2}$ satisfy the following \vspace*{1mm}\\
(i)$ \ \displaystyle l_{n_1,n_2}={-\frac{1}{2} \choose n_1}{-\frac{1}{2} \choose n_2}\frac 1 {\sqrt{-\kappa(0)}}$, \vspace*{3mm} \\
(ii)$\ l_{n_1,n_2}=l_{n_2,n_1}$,\\
where
$${\tau \choose j}=\frac{\tau(\tau-1)\cdots (\tau-j+1)}{j\,!}$$
is the binomial coefficient. }
\end{theorem}
\textbf{Proof.} We write $\lambda(z)=e^{u(z)}$, then $\pa \lambda (z)=\lambda(z)\, \pa u(z)$, and
$$\pa^{n}\lambda(z)=\sum_{j=0}^{{n}-1}{{n}-1 \choose j}\pa^{{n}-j}u(z)
\,\pa^j\lambda(z)$$
by induction, where $\pa^0\lambda(z)=\bar{\pa}^0\lambda(z)=\lambda(z),$
so
$$l_{0,{n_2}}=\frac{1}{{n_2}!}\lim_{z\rightarrow0}\sum_{j=0}^{{n_2}-1}{{n_2}-1 \choose j}
z^{{n_2}-j}\pa^{{n_2}-j}u(z)\cdot|z|\log(1/|z|)z^j\pa^j\lambda(z).$$
Theorem \ref{estimate for cusps original} gives that $l_{0,0}=1/\sqrt{-\kappa(0)}$. From the existence of $\lim_{z\rightarrow0}z^{{n_2}-j}\pa^{{n_2}-j}u(z)$ and $l_{0,0}$, we know that $l_{0,\;{n_2}}$ exists. Next, limit (ii) in Theorem \ref{general u} enables us to write $l_{{n_1},{n_2}}$ as a sum of the terms only containing pure derivatives of $u(z)$,
\be\label{ind}
l_{{n_1},{n_2}}=\frac{1}{{n_1}!{n_2}!}\lim_{z\rightarrow0}\sum_{j=0}^{{n_2}-1}{{n_2}-1 \choose j}
z^{{n_2}-j}\pa^{{n_2}-j}u(z)\,|z|\log(1/|z|)\bar{z}^{n_1} z^j\bar{\pa}^{n_1}\pa^j\lambda(z),
\ee
thus the existence of $l_{0,{n_2}}$ guarantees $l_{{n_1},{n_2}}$ exists.

By Theorem \ref{coro}, it is known that $l_{0,1}$ is a real number, so $l_{1,0}=\overline{l_{0,1}}=l_{0,1}$.
Since
\be \label{pa bar n lam}
{\displaystyle\bar{\pa}^{n_2}\lambda(z)=\sum_{j=0}^{{n_2}-1}{{n_2}-1 \choose j}
\bar{\pa}^{{n_2}-j}u(z)\,\bar{\pa}^j\lambda(z)},
\ee
then $l_{{n_2},0}=l_{0,{n_2}}$ by induction. From \eqref{ind}, \eqref{pa bar n lam}, and (i) of Theorem \ref{general u}, we have
\ben
&&l_{{n_1},{n_2}}\\
&=&\sum_{j=0}^{{n_2}-1} \lim_{z\rightarrow0} \frac{1}{{n_1}!{n_2}!} \frac{({n_2}-1)!}{j!({n_2}-1-j)!}
z^{{n_2}-j}\pa^{{n_2}-j}u(z)\cdot|z|\log(1/|z|)\bar{z}^{n_1} z^j\bar{\pa}^{n_1}\pa^j\lambda(z) \\
&=&\frac{1}{{n_2}}\sum_{j=0}^{{n_2}-1}  \frac{1}{{n_1}!} \frac{1}{j!({n_2}-1-j)!}\lim_{z\rightarrow0}
z^{{n_2}-j}\pa^{{n_2}-j}u(z)\cdot\lim_{z\rightarrow0}|z|\log(1/|z|)\bar{z}^{n_1}z^j\bar{\pa}^{n_1}\pa^j\lambda(z) \\
&=&\frac{1}{{n_2}}\sum_{j=0}^{{n_2}-1} \frac{(-1)^{{n_2}-j}}{2}\frac{1}{{n_1}!j!} \lim_{z\rightarrow0}|z|\log(1/|z|)\bar{z}^{n_1} z^j\bar{\pa}^{n_1}\pa^j\lambda(z)
=\frac{1}{2{n_2}}\sum_{j=1}^{{n_2}-1}(-1)^{{n_2}-j}l_{{n_1},j}.
\een
Then
$$ {n_2} \cdot l_{{n_1},{n_2}}=\frac{1}2\sum_{j=0}^{{n_2}-2}(-1)^{{n_2}-j}\,l_{{n_1},j}-\frac{1}{2}\,l_{{n_1},{n_2}-1}
=-({n_2}-1)l_{{n_1},{n_2}-1}-\frac{1}{2}\,l_{{n_1},{n_2}-1}.$$
Since $l_{0,{n_2}}=l_{{n_2},0}$,
\ben
l_{{n_1},{n_2}}&=&\frac{-\frac{1}{2}-{n_2}+1}{{n_2}}l_{{n_1},{n_2}-1}=
{-\frac{1}{2} \choose {n_2}}l_{{n_1},0}\\
&=&{-\frac{1}{2} \choose {n_2}}l_{0,{n_1}}={-\frac{1}{2} \choose {n_2}}{-\frac{1}{2} \choose {n_1}}l_{0,0}.
\een
Thus (i) is valid and (ii) follows form (i). \hfill $\Box$\\

However, when the order $\alpha<1$, the analogous limit
\be \label{lza}
l':=\lim_{z\rightarrow0}|z|^{\alpha}\lambda(z)
\ee
may also exist but cannot be described only in terms of the curvature of $\lambda(z)$. To discuss the limit \eqref{lza} for an SK-metric, we consider the limit
\be \label{upper lza}
l:=\limsup_{z\rightarrow0}|z|^{\alpha}\lambda(z)
\ee
instead, since the definition \eqref{singularity} of a corner guarantees that $l<\infty$ for $l$ defined above. Based on Theorem \ref{Ahlfors lemma} and Corollary 4.4 in \cite{Rothhyper}, we have the following result corresponding to Theorem \ref{higher lambda for cusps} in the thrice-punctured sphere.
\begin{theorem} \label{limitfor3}
\textsl{Let $0<\alpha,\ \beta <1$ and $0<\gamma \leq 1$ such that $\alpha+\beta+\gamma>2$ and $\lambda(z)$ be an SK-metric on the thrice-punctured Riemann sphere $\widehat{\mathbb{C}}\backslash \{0,\,1,\,\infty\}$ of orders $\alpha,\, \beta,\,\gamma$ at $0$, $1$, $\infty$, respectively, with the curvature $\kappa(z)$. Then the upper limit  \eqref{upper lza} satisfies
\be \label{3lgene}
l \leq \frac{\delta}{1-\delta^2}(1-\alpha),
\ee
where
\ben
\delta=\frac{\Gamma( c)}{\Gamma(2- c)}
\left(\frac{\Gamma(1- a)\Gamma(1- b)\Gamma( a+1- c)
\Gamma( b+1- c)}{\Gamma( a)\Gamma( b)\Gamma( c- a)
\Gamma( c- b)} \right)^{1/2}
\een
with
\ben
a=\frac{\alpha+\beta-\gamma}{2},\ b=\frac{\alpha+\beta+\gamma-2}{2},\  c=\alpha.
\een}
\end{theorem}

Now we consider the upper limit $l$ in the once-punctured unit disk $\mathbb{D}^*$. The following result is evident if we combine Theorem \ref{Ahlfors lemma} with Theorem \ref{maximal}.
\begin{theorem}
\textsl{If $\lambda(z)|dz|$ is an an SK-metric on $\mathbb{D}^*$ with the order $\alpha \in (0,\,1)$ at the origin, then the upper limit $l$ defined in \eqref{upper lza} satisfies $l \leq 1-\alpha$.}
\end{theorem}

When the SK-metric satisfies some stronger continuity assumption, the upper limit $l$ in \eqref{upper lza} will become the limit $l'$ in \eqref{lza} at the origin, which enables us to consider the derivatives of $l$ locally. For instance, if the $\kappa(z)$ and $u(z)$ satisfy the assumption of Theorem \ref{general u}, then $u(z)$ is of class $C^{\;\!2}$ in a neighborhood of $z=0$ and $l=l'$ locally holds near the origin. Therefore we state the following result in terms of the recurrence relation similar to Theorem \ref{higher lambda for cusps}.
\begin{theorem}
\textsl{Let the functions $\kappa(z)$ and $u(z)$ satisfy the assumption of Theorem \ref{general u} for an integer $n$ and let $\lambda(z):=e^{u(z)}$. If the order of $u(z)$ is $\alpha \in (-\infty, \,1)$, then for $n_1$, $n_2\geq 0$, $n_1+n_2 \leq n$, the limit
$$l_{n_1,n_2}:=\frac{1}{n_1!n_2!}\lim_{z \rightarrow0}|z|{\bar{z}}^{n_1} z^{n_2} {\bar{\pa}}^{n_1} \pa^{n_2} \lambda(z)$$
exists and satisfies the following \vspace*{1mm}\\
(i)$ \ \displaystyle l_{n_1,n_2}={-\frac{\alpha}{2} \choose n_1}{-\frac{\alpha}{2} \choose n_2}l'$, \vspace*{3mm} \\
(ii)$\ l_{n_1,n_2}=l_{n_2,n_1}$, \vspace*{1mm} \\
where $l'$ is defined by \eqref{lza}.}
\end{theorem}

\vspace*{5mm}
\hspace*{-17pt}\textbf{Acknowledgement.} I would like to thank Prof. Toshiyoki Sugawa for his helpful comments, suggestion and encouragement. I also want to thank Prof. Toshihiro Nakanishi for his suggestion on the SK-metrics.\\

\providecommand{\bysame}{\leavevmode\hbox to3em{\hrulefill}\thinspace}
\providecommand{\MR}{\relax\ifhmode\unskip\space\fi MR }
\providecommand{\MRhref}[2]{%
  \href{http://www.ams.org/mathscinet-getitem?mr=#1}{#2}
}
\providecommand{\href}[2]{#2}

\newpage
\end{document}